\documentclass{article}
\UseRawInputEncoding
\pdfoutput=1
\usepackage{amsmath}
\usepackage{amsthm}
\usepackage{amssymb}
\usepackage{color}
\usepackage{mathrsfs}
\usepackage{graphicx}
\usepackage{amsmath}
\numberwithin{equation}{section}
\def\pl{\partial}

\def\*{\raisebox{.5mm}{*}}

\def\t#1{{\widetilde#1}}

\def\dis{\displaystyle}

\def\R{\mathbb{R}}

\def\H{{\cal H}}

\def\PP{{\cal P}}

\def\|{{\Big|}}
\def\({{\Big(}}
\def\){{\Big)}}
\def\[{{\Big[}}
\def\]{{\Big]}}

\def\Ga{\Gamma}

\def\Del{\Delta}

\def\ms{\medskip}

\def\<{\langle}
\def\>{\rangle}

\def\Om{\Omega}

\def\var{\varphi}

\def\det{\mbox{det\,}}

\def\beq{\arraycolsep=1.5pt\begin{eqnarray}}
	\def\eeq{\end{eqnarray}}
\def\bma{\left[\begin{array}}
	\def\ema{\end{array}\right]}
\def\bda{\left|\begin{array}}
	\def\eda{\end{array}\right|}
\def\be{\begin{equation}}
	\def\ee{\end{equation}}

\newtheorem{thm}{{}\hskip\parindent Theorem}[section]
\newtheorem{lem}{{}\hskip\parindent Lemma}[section]

\newtheorem{rem}{{}\hskip\parindent Remark}[section]
\begin{document}
	\large
	\title{Inverse problem of recovering a time-dependent nonlinearity appearing in third-order
		nonlinear acoustic equations\thanks{This work is supported by the National Nature Science Foundation of China under  Grant 12071463, and Postdoctoral Science Foundation of China under Grants 2021TQ0353 and 2022M713315.} }
	
	
	\author{
		Song-Ren Fu\thanks{ Key Laboratory of Systems and Control, Institute of Systems Science, Academy of Mathematics and Systems Science, Beijing 100049, People's Republic of China. 
},  \
		Peng-Fei Yao\thanks{School of Mathematical Sciences, Shanxi University, Taiyuan 030006, People's Republic of China.  }  \ and\
	 Yongyi Yu\thanks{ Corresponding  author.  School of Mathematics, Jilin University, Changchun 130012, People's Republic of China.  E-mail address: yuyy122@amss.ac.cn.}}
	
	\date{}
	
	\maketitle
	
	{\bf Abstract.} This paper is devoted to some inverse problems of recovering the nonlinearity for a Jordan-Moore-Gibson-Thompson equation (J-M-G-T equation for short), which is a third order nonlinear acoustic equation.
The well-posedness of the nonlinear equation is obtained for the small initial and boundary data. By the second order linearization to the nonlinear equation, and construction of complex geometric optics (CGO for short) solutions for the linearized equation, the uniqueness of recovering the nonlinearity is derived. 

\ms
	

	\section{Introduction}
Nonlinear wave propagation in acoustic environment has  been an
interesting and active topic.
It is widely used in medical and industrial applications,  including medical imaging,  thermotherapy, ultrasound cleaning and sonochemistry, see \cite{OVA, MHJR, ku}.
The classical models of nonlinear acoustics are second order in time, and known as the Kuznetsov's equation, the Westervelt equation, and the Kokhlov-Zabolotskaya-Kuznetsov equation.

As we all known, by the standard
Fourier temperature law, a thermal disturbance at one point has an instantaneous
effect elsewhere in the medium, which may leads to an
infinite signal speed paradox.
 In the context of nonlinear acoustic waves
modeling high-intensity ultrasound, replacement of the Fourier law by the Maxwell-Cattaneo law yields the J-M-G-T equation, which is a third order (in time) nonlinear equation.
The linear version is called the Moore-Gibson-Thompson equation (M-G-T equation for short).
This paper is devoted to investigating the inverse problem of determining the time-dependent nonlinearity term for the J-M-G-T equations.

	Let $\Om\subset \R^n$ be a bounded domain of dimension $n\ge 2$ with smooth boundary $\pl\Om=\Ga.$ For given $T>0,$ let $Q=\Om\times (0,T)$ and $\Sigma=\Ga\times(0,T).$
	We consider the following nonlinear J-M-G-T equation
	\begin{eqnarray}\label{sysmgt}
		\left\{ \begin{array}{ll}
			\tau u_{ttt}+\alpha u_{tt}-c^2\Del u-b\Del u_t=F(u,u_t,u_{tt},\nabla u,\nabla u_t), \ & (x,t)\in Q,\\
			u=h(x,t), \ & (x,t)\in \Sigma,\\
			u(x,0)=u_0(x), u_t(x,0)=u_1(x), u_{tt}(x,0)=u_2(x),\  & x\in\Om,
		\end{array} \right.
	\end{eqnarray}
	where $\tau, c, \alpha, b$ refer to the thermal relaxation coefficient, the speed, the viscosity parameter (friction damping) and diffusivity of sound, respectively.
This third order in time equation arises, for example,  from the wave propagation in viscous thermally relaxing fluids.
For more details of the above model, we refer to \cite{PMJ, mg, T} for its derivation.
We also refer to \cite{KLP, LZ,VB,RRB} for the nonlinear acoustic models on well-posedness, stability and controllability.

	
There are two kinds of nonlinearities for J-M-G-T equation (\ref{sysmgt}).
	The first case is Westervelt-type, where
	\be\label{first}
	F(u,u_t,u_{tt},\nabla u,\nabla u_t)=(p(x,t)u^2)_{tt}=p_{tt}u^2+4p_tuu_t+2p(u_t^2+uu_{tt}).
	\ee
	The second is
	\be\label{second}
	F(u,u_t,u_{tt},\nabla u,\nabla u_t)=(pu_t^2+k|\nabla u|^2)_t,
	\ee
	which refers to the Kuznetsov-type. In both cases, $p(x,t)$ is the time-dependent parameter of nonlinearity that arises from e.g., the pressure-density relation in a given medium.
	
	In this paper, we study the inverse problem of recovering the nonlinearity for Westervelt-type J-M-G-T equation. The time-dependent parameter of nonlinearity $p(x,t)$ is uniquely recovered from knowledge of the Dirichlet-to-Neumann  (DtN for short) map. While for the second case (\ref{second}), we leave it as our future work.

	The inverse problem of determining physical parameters of PDEs have attracted much attention with lots of literature on inverse elliptic equations, parabolic equations, Schr\"odinger equations, hyperbolic equations, etc. Among those, inverse problems for nonlinear acoustics have raised recently. In \cite{GUYZ}, the nonlinear coefficient is recovered from the DtN map for the Westervelt equation without damping. Later, the inverse problem for a nonlinear acoustic equation with a damping term and a general nonlinearity was recently studied in \cite{ZY}. The author proved that the one-form
	and the nonlinearity can be determined from the DtN map, up to a gauge transformation. We also refer to \cite{SGJ}, in which the authors concerned the ultrasound imaging problem governed by a nonlinear wave equation of Westervelt-type with variable wave speed. 

In this paper,
a second order linearization and the CGO solutions are used to derive the main result. It is pointed out that, the higher order linearization method, proposed in \cite{VI}, is very useful in dealing with the nonlinear equations.  The CGO solutions are also crucial in inverse problems. We refer to e.g., \cite{AFYK, AJY, PGJ, LLL, GUJZ} for the using of these methods. Other inverse problems of nonlinear acoustics can be referred to \cite{NEPS, AJ, BKVN, LLZ} and references therein.
	
	For the inverse problem of J-M-G-T equation, even for the linear M-G-T equation, few results are available in literature. 
In \cite{Kal}, the authors proved the uniqueness of the space-dependent nonlinearity parameter for J-M-G-T equation by the Inverse Function Theorem.  The observation is a single time trace of the acoustic pressure on the boundary.
	In \cite{RRAS,LiuT}, the authors studied the inverse M-G-T equation of determining the space varying frictional damping term $\alpha(x)$ by the Carleman estimate.
We also mention that,  a stability result of recovering the time-dependent coefficient $\alpha(x,t)$ is obtained in \cite{fy} by means of the construction of geometric optics solutions.
	To the authors' best knowledge,  the inverse problem of determining the time-dependent nonlinearity term in the nonlinear J-M-G-T equation has not been studied yet.

The rest of this paper is organized as follows: In Section 2, we state the main results for the inverse problems. In Section 3, the well-posedness of the nonlinear J-M-G-T equation is proved. Section 4 concerns to construct the CGO solutions for the M-G-T equation. In Section 5, we prove the main theorems of this present paper. Finally, Section 6 gives the concluding remarks.

\section{Statement of main results}	
	In this section, we give some preliminaries and the main results.
	For simplicity, we assume that $\tau$ is normalized to be $1$, the parameters $\alpha, b, c$ are positive constants, and also assume that there exists a positive constant $M>1$ such that
	\be\label{MMM}
	M^{-1}\le \alpha, b, c\le M.
	\ee
	Denote by $\gamma=\alpha-\frac{c^2}{b}.$ It is well known that $\gamma$ is a crucial parameter that related to the exponential stability, energy conservation and unstable of the system.
	
	Here, some function spaces are introduced.
	For any non-negative integer $m,$ let
	\be\label{Hm}
	H^m(\Sigma)=\mathop\cap_{k=0}^mH^k(0,T;H^{m-k}(\Gamma)),\quad\H_{m}=H^{m+2}(\Omega)\times H^{m+1}(\Omega)\times H^m(\Omega),
	\ee
	\be\label{boundarydata}
	A_m=C^m([0,T];H^{m+\frac32}(\Gamma))\cap C^{m+1}([0,T];H^{m+\frac12}(\Gamma))\cap H^{m+2}(\Sigma),
	\ee
	\be\label{Nm1}
	N_{m+1}=H^{m+1}(\Sigma)\backslash L^2(0,T;H^{m+1}(\Gamma))=\mathop\cap_{k=1}^{m+1}H^k(0,T;H^{m+1-k}(\Gamma)).
	\ee
	Moreover, we denote by
	$$E^m=\mathop\cap_{k=0}^mC^k([0,T];H^{m-k}(\Omega)),$$
	with norm
	$$||u||_{E^m}=\sup_{0\le t\le T}\sum_{k=0}^{m}||\partial_t^ku(\cdot,t)||_{H^{m-k}(\Omega)}.$$
	It is well known that $E^m$ is a Banach algebra provided that $m>n+1,$ where $n$ is the dimension of $\Omega.$

	Let $u(x,t)=u(x,t;h,u_0,u_1,u_2,p)$ denote the solution of (\ref{sysmgt}) with respect to the boundary and initial data $(h,u_0,u_1,u_2)$, and the parameter $p(x,t).$ The DtN map is defined as
	$$\mathcal D: H^{s+1}(\Sigma)\to H^s(\Sigma),\quad \mathcal D h=\frac{\partial u}{\partial\nu}\Big|_{\Sigma},\quad s\in\mathbb R,$$
	where $\nu$ is the outward unit normal vector filed along $\Gamma.$

First, we present the first inverse problem on determining the time-dependent parameter of nonlinearity by all boundary and initial values. For given $T>0,$ define the following input-to-output map by
	$$\Lambda_T: A_m\times\H_{m}\to N_{m+1}\times \H_m,$$
	$$\Lambda_T(h,u_0,u_1,u_2)=(\mathcal Dh, u(x,T), u_t(x,T),u_{tt}(x,T)).$$

In this paper, we are first interested in the following inverse problem:

\ms

	{\bf $\bullet$ Inverse Problem 1.} Let $(h,u_0,u_1,u_2)$ be system inputs. Recovering $p(x,t)$ by $\Lambda_T.$
	
\ms

Here, some assumption conditions are introduced for the above inverse problem.
	
\ms

	{\bf (A.1)} Let $m>n+1.$ The parameter $p\in E^{m+2}$ such that
	$$||p||_{E^{m+2}}\le M,$$
	where the positive constant $M$ is the same with that in (\ref{MMM}). Moreover, we assume that
there exists a small constant $\delta >0$ such that	
\be\label{c1}
(h,u_0,u_1,u_2)\in A_m\times \H_m
\ee
	satisfying
\be\label{c2}
||u_0||_{H^{m+2}(\Omega)}+||u_1||_{H^{m+1}(\Omega)}+||u_2||_{H^m(\Omega)}+||h||_{A_m}\le \delta.
\ee
For a given sufficiently small $\t\delta>0$, let
	\begin{equation*}\label{obserT}
	T^*:={\rm diam}\ \Omega+\t\delta,
	\end{equation*}
	where ${\rm diam}\ \Omega$ denotes the diameter of $\Omega.$
We choose a domain $\Omega_1$ such that $\Omega\subset\subset\Omega_1$ satisfying
$${\rm diam}\ \Omega<{\rm diam}\ \Omega_1\le {\rm diam}\ \Omega+\t\delta.$$
For more details on $\Om_1$, see Section 4.

\ms

	We are now in a position to state the first result.
	\begin{thm}\label{thm1}
		Let assumption {\bf (A.1)} hold. Then, for given $T>0,$ the map $\Lambda_T$ uniquely determines $p(x,t)$ for $(x,t)\in \Omega\times(0,T).$
	\end{thm}

\begin{rem}
From the above result, it is clear that the parameter $p(x,t)$ can be uniquely determined by $\Lambda_T$ in the full domain $Q$ with all possible boundary and initial data as the system inputs.
\end{rem}
More importantly, we also want to determine $p(x,t)$ only by the boundary data for the nonlinear equation (\ref{sysmgt}). Similarly, another map is given via
	$$\mathcal B_T: A_m\to N_{m+1}\times \H_m,$$
	$$\mathcal B_T(h)=(\mathcal Dh, u(x,T), u_t(x,T),u_{tt}(x,T)).$$

From the above definitions, $\Lambda_T$ acts on all boundary and initial data, while, $\mathcal B_T$ acts on only the boundary data $h.$
The maps $\Lambda_T$ and $\mathcal B_T$ are well-defined and continuous (see Theorem \ref{wellnonlinear} below). We are concerned with the following inverse problem.

\ms

	{\bf $\bullet$ Inverse Problem 2.} Let $h$ be the boundary input. Recovering $p(x,t)$ by $\mathcal B_T.$
	
\ms

	The following assumptions is given.
	
	\ms

{\bf (A.2)} Let $m>n+1.$ The parameter $p\in E^{m+2}$ such that
	$$||p||_{E^{m+2}}\le M,\quad {\rm supp}\ p\subset\overline\Omega\times(T^*+2\t\delta,T).$$
 We also assume that $(\ref{c1})$ and $(\ref{c2})$ hold here.
	
	\ms
	
	The main result  established for {\bf Inverse Problem 2} is stated as follows:
	\begin{thm}\label{thm2}
		Let assumption {\bf (A.2)} hold and $T>T^*+2\t\delta.$ Then, the map $\mathcal B_T$ uniquely determines
		$p(x,t)$ for $(x,t)\in \Omega\times(0,T).$
	\end{thm}
	\begin{rem}
		{ 
From Theorem {\rm\ref{thm2}} and the assumption {\bf (A.2)}, we know that $p(x,t)$ can be determined in $\Omega\times(T^*+2\t\delta,T)$ only by the Dirichlet boundary data. The main reason is  the cut-off procedure used in the proof.

It is pointed out that, as we will see in the subsequent discussions, the uniqueness results in Theorems 1.1 and 1.2 are also valid for variable coefficient $b(x)$ provided that $g=b^{-1}ds^2$ is a simple Riemannian metric.
		}
	\end{rem}

\section{Well-posedness of the nonlinear system (\ref{sysmgt})}
This section is devoted to giving the well-posedness of the nonlinear J-M-G-T equation (\ref{sysmgt}) with small initial and boundary data.
We begin with the following linear M-G-T equation:
	\begin{equation}\label{linear}
		\left\{\begin{array}{ll}
			v_{ttt}+\alpha v_{tt}-c^2\Del v-b\Del v_t=f(x,t), \ & (x,t)\in Q,\\
			v=\xi(x,t),\  & (x,t)\in \Sigma,\\
			v(x,0)=v_0(x), v_t(x,0)=v_1(x), v_{tt}(x,0)=v_2(x),\ &  x\in \Omega.
		\end{array}
		\right.
	\end{equation}
	Here and after, we briefly denote $L=\partial_t^2-b\Delta$ and
	$$\PP =\partial_t^3+\alpha\partial_t^2-b\Delta\partial_t-c^2\Delta.$$
	Then,
	\begin{equation*}\label{Lv}
	\PP v=Lv_t+\beta Lv+\gamma v_{tt},\quad (x,t)\in Q,
	\end{equation*}
	where $\beta=\frac{c^2}{b}>0$ is a constant.
	
	\ms

	Let $\left\langle\cdot,\cdot\right\rangle$ be the standard Euclidean inner product throughout this paper.
	We will use many times the simple Riemannian metric $g=b^{-1}ds^2.$  Let $\nabla,\nabla_g$ respectively be the gradient operator in the Euclidean metric and the Riemannian metric $g.$ Then for any function $\psi, \psi_1,\psi_2\in H^1(\Omega),$ we have
	$$\nabla_g\psi=b\nabla \psi,\quad \left\langle \nabla_g\psi_1,\nabla_g\psi_2\right\rangle _g=g(\nabla_g\psi_1,\nabla_g\psi_2)=\left\langle b\nabla\psi_1,\nabla\psi_2\right\rangle,$$
	where  $\left\langle\cdot,\cdot \right\rangle_g=g(\cdot,\cdot)$ denotes the inner product in the metric $g.$

	We denote by $D$ the Levi-Civita connection in the metric $g$,  $DH$ is the  covariant derivative of a certain vector filed $H,$ which is a tensor field of second order. Besides, denote by ${\rm div}$ the standard divergence operator in the Euclidean metric. It is pointed out that, though the Riemannian metric $g$ here is simple,  the arguments can be conveniently carried in the Riemannian setting.
	
	\ms
	
	Set $$E_0(v(t))=\frac12\int_\Omega(v_t^2+|\nabla v|^2)dx,\ t\ge 0.$$
	Define the energy associated to system (\ref{linear}) as
	\begin{eqnarray*}\label{Nm}
	E(t)&&=E_0(v(t))+E_0(v_t(t))+\frac12\int_\Omega|\Delta v|^2dx\nonumber\\
	&&=\frac12\int_\Omega(v_t^2+v_{tt}^2+|\nabla v|^2+|\nabla v_t|^2+|\Delta v|^2)dx,\ t\ge 0.
	\end{eqnarray*}
	For the linear system (\ref{linear}), from \cite[Theorem 1.1]{FBME}, we have the following lemma:
	\begin{lem}\label{well-lower}
		Let $T\in (0,+\infty)$ be given. Assume that
		$$(v_0,v_1,v_2)\in H^2(\Omega)\times H^1(\Omega)\times L^2(\Omega),\quad f\in L^2(Q),$$
		and
		$$\xi\in C([0,T];H^{\frac32}(\Gamma))\cap H^2(0,T;L^2(\Gamma)),\quad \xi_t\in C([0,T];H^{\frac12}(\Gamma)),$$
		satisfying the compatibility conditions:
		$$v_0|_\Gamma=\xi(x,0)\in H^{\frac32}(\Gamma),\quad v_1|_\Gamma=\xi_t(x,0)\in H^{\frac12}(\Gamma).$$
		Then the linear system {\rm(\ref{linear})} admits a unique solution
		$$(v,v_t,v_{tt})\in C([0,T];H^2(\Omega)\times H^1(\Omega)\times L^2(\Omega)).$$
		If in addition, $\xi\in H^2(\Sigma),$ then we have the hidden regularity
		$$\frac{\partial^2v}{\partial t\partial\nu}\in L^2(\Sigma).$$
		
	\end{lem}
	
	In order to get a higher regularity of the linear system (\ref{linear}), which will be used in the analysis of the inverse problems, we begin with the compatibility condition of order $m+1.$ More precisely, we need the following
	\begin{equation}\label{compatibility}
		\left\{\begin{array}{l}
			\partial_t^kv(x,0)=\partial_t^k\xi(x,0)\quad {\rm on}\ \Gamma,\quad 0\le k\le 2,\vspace{1.0ex} \\
			\partial_t^{k}\xi(x,0)=-\alpha\partial_t^{k-1}v(x,0)+b\Delta\partial_t^{k-2}v(x,0)\\
			\quad \qquad \qquad +c^2\Delta\partial_t^{k-3}v(x,0)+\partial_t^{k-3}f(x,0)\quad {\rm on}\ \Gamma,\quad 3\le k\le m+1.
		\end{array}
		\right.
	\end{equation}

	\begin{thm}\label{well-linear}
		Let $m\ge0$ be an integer. Suppose that $f\in E^m$ and
		$$(v_0,v_1,v_2,\xi)\in \H\times A_m,$$
		satisfying the compatibility conditions {\rm(\ref{compatibility})} up to order $m+1.$ Then,
		the linear system {\rm(\ref{linear})} admits a unique solution $v\in E^{m+2}.$ Moreover, there exists a positive constant $C=C(\Omega, T,M)$ such that
		\beq\label{estlinear}
		&&||v||_{E^{m+2}}+||\partial_\nu v||_{N_{m+1}}\nonumber\\
		&&\le C\left(||v_0||_{H^{m+2}(\Omega)}+||v_1||_{H^{m+1}(\Omega)}+||v_2||_{H^{m}(\Omega)}+||\xi||_{H^{m+2}(\Sigma)}+||f||_{E^m}\right).
		\eeq
	\end{thm}
\begin{rem}
The above result is based on Lemma {\rm\ref{well-lower}}.
It's also an interesting problem whether the well-posedness can be achieved provided that the Dirichlet data $\xi\in H^{m+2}(\Sigma).$ For the classical wave equation of second order, the Dirichlet data is assumed to be $H^{m+1}(\Sigma)$ in obtaining the $E^{m+1}$ regularity of solution, e.g., see {\rm\cite[Theorem 2.45]{AKYM}.
}
\end{rem}
	{\bf Proof of Theorem \ref{well-linear}}\ The proof is divided into three steps.
	
	{\bfseries Step 1.}\ We begin with the case where $m=0.$ By Lemma \ref{well-lower},  the linear system {\rm(\ref{linear})} admits a unique solution
	$$v\in C([0,T];H^2(\Omega))\cap C^1([0,T];H^1(\Omega))\cap C^2([0,T];L^2(\Omega)).$$
	
	Respectively multiplying both sides the first equation of (\ref{linear}) by $v_{tt}$ and $v_t$,  integrating over $\Omega\times(0,t),$ applying the Green formula, we have
	\begin{eqnarray*}
	&&\frac12\int_\Om (v_{tt}^2+b|\nabla v_t|^2)dx+\alpha\int_0^t\int_\Om v_{tt}^2dxdl+c^2\int_\Om\left\langle\nabla v,\nabla v_t \right\rangle dx\nonumber\\
	&&=\frac12\int_\Om(v_2^2+b|\nabla v_1|^2)dx+c^2\int_\Om\left\langle\nabla v_1,\nabla v_0\right\rangle dx+c^2\int_0^t\int_\Om|\nabla v_t|^2dxdl\nonumber\\
	&&\quad +\int_0^t\int_\Om fv_{tt}dxdl+\int_0^t\int_\Ga \xi_{tt}(c^2\partial_\nu v+b\pl_\nu v_t)]d\Ga dl,
	\end{eqnarray*}
	and
	\begin{eqnarray*}
	&&\frac12\int_\Om(\alpha v_t^2+c^2|\nabla v|^2)dx+b\int_0^t\int_\Om |\nabla v_t|^2dxdl\nonumber\\
	&&=\int_0^t\int_\Om v_{tt}^2dxdl+\frac12\int_\Om(c^2|\nabla v_1|^2+2v_1v_2)dx\nonumber\\
	&&\quad +\int_0^t\int_\Om fv_tdxdl+\int_0^t\int_\Ga \xi_t(c^2\pl_\nu v+b\pl_\nu v_t)d\Ga dl.\nonumber
	\end{eqnarray*}
	Add the above two equalities, and notice that $\alpha,b,c$ are  constants with positive lower bound. Then, by the Cauchy-Schwartz inequality, we have
	\beq\label{est0}
	&&\frac12\int_\Om(v_{tt}^2+b|\nabla v_t|^2+\alpha v_t^2+c^2|\nabla v|^2)dx\nonumber\\
	&&\le C\left( ||v_0||^2_{H^1(\Om)}+||v_1||^2_{H^1(\Om)}+||v_2||^2_{L^2(\Om)}\right) +C\int_0^t\int_\Om(|\nabla v_t|^2+v_{tt}^2+v_t^2) dxdl\nonumber\\
	&&\quad +\frac{b}{4}\int_\Om |\nabla v_t|^2dx+\frac{4c^4}{b}\int_\Om|\nabla v|^2dx\nonumber\\
	&&\quad +C\left( ||f||^2_{L^2(Q)}+||\xi||^2_{H^2(0,T;L^2(\Ga))}+|||\pl_\nu v||^2_{H^1(0,T;L^2(\Ga))}\right).
	\eeq
	
	Notice that
	\be
	\nabla v(x,t)=\nabla v(x,0)+\int_0^t\nabla v_t(x,l)dl.\nonumber
	\ee
	Then, there exists a positive constant $C=C(\Omega,T)$ such that
	\be\label{calcus}
	\int_\Om|\nabla v|^2dx\le C||v_0||^2_{H^1(\Omega)}+C\int_0^t\int_\Om|\nabla v_t|^2dxdl.
	\ee
	
	Multiplying the equation (\ref{linear}) by $Lv$ and integrating over $\Omega\times(0,t),$ we have
	\beq\label{est1}
	&&\frac12\int_\Om|Lv|^2dx+\frac\beta 2\int_0^t\int_\Om|Lv|^2dxdl\nonumber\\
	&&\le \frac{4\gamma^2}{\beta}\int_0^t\int_\Om v_{tt}^2dxdl+\frac4 \beta ||f||^2_{L^2(Q)}+||v_2||^2_{L^2(\Om)}+b^2||\Del v_0||^2_{L^2(\Om)}.
	\eeq
	Together with (\ref{est0})-(\ref{est1}) implies that
	\beq\label{est2}
	&&\frac12\int_\Om(v_{tt}^2+|\nabla v_t|^2+v_t^2+|\nabla v|^2+|\Del v|^2)dx\nonumber\\
	&&\le C\left( ||v_0||^2_{H^2(\Om)}+||v_1||^2_{H^1(\Om)}+||v_2||^2_{L^2(\Om)}\right) +C\int_0^t\int_\Om(|\nabla v_t|^2+v_{tt}^2+v_t^2) dxdl\nonumber\\
	&&\quad +C\left( ||f||^2_{L^2(Q)}+||\xi||^2_{H^2(0,T;L^2(\Ga))}+||\pl_\nu v||^2_{H^1(0,T;L^2(\Ga))}\right).
	\eeq
	
	{\bfseries Step 2.}\ We next to estimate the term $||\pl_\nu v||^2_{H^1(0,T;L^2(\Ga))}$ appearing in the last term of (\ref{est2}).
	
	Recalling the classical multipliers $Hv=\left\langle H,\nabla v\right\rangle $ and $Hv_t=\left\langle H,\nabla v_t \right\rangle $, where
	$$H=\sum_{k=1}^{n}H_k(x)\frac{\partial}{\partial x_k},\ H_k(x)\in C^1(\overline\Omega),\ \forall k=1,2,\cdots,n,\ x=(x_1,x_2,\cdots,x_n)\in\mathbb R^n.$$
	That is, $H$ is a smooth vector field on $\overline{\Omega}.$
	
	We calculate $2(Hv)Lv=2(Hv)(v_{tt}-b\Delta v)$ and integrate it over $\Omega\times(0,t).$
	
	Indeed, by \cite[Theorem 2.1]{Yaobook}, taking $A(x)=bI_n,$ with $I_n$ the $n$-th identity matrix, we have
	\beq\label{multiplier1}
	2(Hv)Lv&&=2(v_tHv)_t+{\rm div}[(|\nabla_g v|_g^2-v_t^2)H-2(Hv)\nabla_g v]\nonumber\\
	&&\quad +2DH(\nabla_g v,\nabla_g v)+(v_t^2-|\nabla_g v|_g^2){\rm div}H,
	\eeq
	where $g=b^{-1}ds^2,$ the Riemannian metric.
	
	Integrating (\ref{multiplier1}) over $\Omega\times(0,t),$ via the divergence theorem, we can get
	\beq\label{mulint}
	2\int_0^t\int_\Om(Hv)Lvdxdl&&=2\int_\Om v_tHvdx-2\int_\Om v_1Hv_0dx\nonumber\\
	&&\quad+\int_0^t\int_\Om[2DH(\nabla_g v,\nabla_g v)+(v_t^2-|\nabla_g v|_g^2){\rm div}H]dxdl\nonumber\\
	&&\quad +\int_0^t\int_\Ga[(b|\nabla v|^2-\xi_t^2)\left\langle H,\nu\right\rangle-2b(Hv)\pl_\nu v]d\Ga dl,
	\eeq
	where $\tau$ is the unit tangent vector field along $\Gamma$ with respect to the standard Euclidean metric $I_n,$ and $|\nabla v|^2=|\tau v|^2+|\partial_\nu v|^2$ on $\Gamma.$
	
	Taking $H=\nu$ on $\Gamma$ in (\ref{mulint}), we obtain
	\beq\label{mul1}
	b\int_0^t\int_\Ga|\pl_\nu v|^2d\Gamma dl&&=-2\int_0^t\int_\Om(Hv)Lvdxdl+2\int_\Om v_tHvdx-2\int_\Om v_1Hv_0dx\nonumber\\
	&&\quad +\int_0^t\int_\Om[2DH(\nabla_g v,\nabla_g v)+(v_t^2-|\nabla_g v|_g^2){\rm div}H]dxdl\nonumber\\
	&&\quad +\int_0^t\int_\Ga (b|\tau \xi|^2-\xi_t^2)d\Ga dl.
	\eeq
	
	Multiplying $Lv_t$ by $Hv_t,$ by the same argument, we have
	\beq\label{mul2}
	b\int_0^t\int_\Ga|\pl_\nu v_t|^2d\Gamma dl&&=-2\int_0^t\int_\Om(Hv_t)Lv_tdxdl+2\int_\Om u_{tt}Hv_tdx-2\int_\Om v_2Hv_1dx\nonumber\\
	&&\quad +\int_0^t\int_\Om[2DH(\nabla_g v_t,\nabla_g v_t)+(v_{tt}^2-|\nabla_g v_t|_g^2){\rm div}H]dxdl\nonumber\\
	&&\quad +\int_0^t\int_\Ga (b|\tau \xi_t|^2-\xi_{tt}^2)d\Ga dl.
	\eeq
	Notice that
	$$Lv_t=f-\gamma v_{tt}-\beta Lv,\ |Hv|\le C|\nabla v|,\ |Hv_t|\le C|\nabla v_t|,$$ and
	$$DH(\nabla_gv,\nabla_gv)\le C|\nabla v|^2,\ DH(\nabla_gv_t,\nabla_gv_t)\le C|\nabla v_t|^2,$$
	where the constant $C$ depends on $b$ and $H.$
	It then follows from (\ref{est2}), (\ref{mul1}) and (\ref{mul2}) that
	\beq\label{energy}
	&&\int_\Om(v_{tt}^2+|\nabla v_t|^2+v_t^2+|\nabla v|^2+|\Del v|^2)dx\le C\left( ||v_0||^2_{H^2(\Om)}+||v_1||^2_{H^1(\Om)}+||v_2||^2_{L^2(\Om)}\right)\nonumber\\
	&&\quad +C\int_0^t\int_\Om(|\nabla v_t|^2+v_{tt}^2+v_t^2+|\nabla v|^2+|\Del v|^2) dxdl\nonumber\\
	&&\quad +C\left( ||f||^2_{L^2(Q)}+||\xi||^2_{H^2(0,T;L^2(\Ga))}+||\xi||^2_{H^1(0,T;H^1(\Ga))}\right).
	\eeq
	Applying the Gronwall inequality and combining with (\ref{mul2})-(\ref{energy}), yields
	\beq
	&&E(t)+||\pl_\nu v||_{H^1(0,T;L^2(\Sigma))}\nonumber\\
	&&\le C\left(||v_0||^2_{H^2(\Om)}+||v_1||^2_{H^1(\Om)}+||v_2||^2_{L^2(\Om)}+||f||^2_{L^2(Q)}+||\xi||^2_{H^2(\Sigma)}\right).\nonumber
	\eeq
	
	{\bfseries Step 3.} Finally, we improve the regularity by the bootstrap arguments.
	Assume that
	$$f\in L^2(0,T;H^m(\Omega)),\ \partial_t^mf\in L^2(0,T;L^2(\Omega)),\ (v_0,v_1,v_2,\xi)\in \H\times A_m.$$
	If moreover, the compatibility condition (\ref{compatibility}) holds, it is clear that system (\ref{linear}) admits a unique solution
	$$v\in C([0,T];H^{m+2}(\Omega))\cap C^{m+2}([0,T];L^2(\Omega)),$$
	and
	$$\partial_\nu v\in Y_{m+1}:=H^{m+1}(0,T;L^2(\Gamma))\cap H^1(0,T;H^m(\Gamma))$$
	satisfying
	\beq
	&&||v||^2_{C([0,T];H^{m+2}(\Om))}+||v||^2_{C^{m+2}([0,T];L^2(\Om))}+||\pl_\nu v||^2_{Y_{m+1}}\nonumber\\
	&&\le C\left(||v_0||^2_{H^{m+2}(\Om)}+||v_1||^2_{H^{m+1}(\Om)}+||v_2||^2_{H^m(\Om)}+||\xi||^2_{A_m}\right.\nonumber\\
	&&\left. \quad\quad\quad+||f||^2_{L^2(0,T;H^m(\Om))}+||\pl_t^mf||^2_{L^2(Q)}\right).
	\eeq
	
	We next to show $v\in E^{m+2}.$ For any integer $k\in [1,m+2],$ it suffices to prove that
	$$\partial_t^kv(\cdot,\cdot)\in C([0,T];H^{m+2-k}(\Omega)).$$	
	Set $\hat v=\partial_t^kv.$  Then $\hat v$ solves
	\be\label{part}
	\left\{\begin{array}{ll}
		\hat v_{ttt}+\alpha \hat v_{tt}-b\Del \hat v_t-c^2\Del \hat v=\partial_t^kf,\  &  (x,t)\in Q,\\
		\hat v=\partial_t^k\xi,\  & (x,t)\in\Sigma.
	\end{array}
	\right.
	\ee
	By the assumption that
	$$f\in E^m,\ (v_0,v_1,v_2,\xi)\in \H_m\times A_m,$$
	we have
	$$\partial_t^kf\in L^2(0,T;H^{m-k}(\Omega)),\quad \partial_t^k\xi\in A_{m-k},$$ $$(\hat v(\cdot,0),\hat v_t(\cdot,0),\hat v_{tt}(\cdot,0))\in H^{m+2-k}(\Omega)\times H^{m+1-k}(\Omega)\times H^{m-k}(\Omega).$$
	Additionally, the compatibility condition for (\ref{part}) of order $m+1-k$ holds. Therefore,
	by the above arguments, we have
	$$\hat v=\partial_t^kv\in C([0,T];H^{m+2-k}(\Omega)),\ k=1,2,\cdots,m+2.$$
	Therefore, we finish the proof of Theorem \ref{well-linear}.  \qedsymbol
	
	\ms
	
	We are now in a position to prove the well-posedness of the nonlinear system (\ref{sysmgt}).
	\begin{thm}\label{wellnonlinear}
		Assume that the integer $m>n+1$, $(u_0, u_1, u_2, h, p)\in \H_m\times A_m\times E^{m+2}$ and $(\ref{c2})$ holds.
		Also, suppose that the compatibility condition {\rm(\ref{compatibility})} holds for {\rm(\ref{sysmgt})}, in which $f(x,t)$ is replaced by $(pu^2)_{tt}$.
		Then, the nonlinear system {\rm(\ref{sysmgt})} admits a unique solution $u\in E^{m+2},$ and there exists a positive constant $C=C(\Omega,T,M)$ such that
		\beq\label{estnon}
		&&||u||_{E^{m+2}}+||\pl_\nu u||_{N_{m+1}}\nonumber\\[2mm]
		&&\le C\left( ||u_0||_{H^{m+2}(\Omega)}+||u_1||_{H^{m+1}(\Omega)}+||u_2||_{H^m(\Omega)}+||h||_{H^{m+2}(\Sigma)}\right).
		\eeq
	\end{thm}
	{\bf Proof.}\quad Let $v_1$ solves
	\begin{equation}\label{sysv1}
		\left\{\begin{array}{ll}
			v_{1ttt}+\alpha v_{1tt}-b\Delta v_{1t}-c^2\Delta v_1=0,\ & (x,t)\in Q,\\
			v_1=\xi,\ & (x,t)\in \Sigma,\\
			(v_1(x,0),v_{1t}(x,0),v_{1tt}(x,0))=(u_0,u_1,u_2),\ & x\in \Omega.
		\end{array}
		\right.
	\end{equation}
	By Theorem \ref{well-linear}, there exists a unique solution $v_1\in E^{m+2}$ such that
	$$||v_1||_{E^{m+2}}+||\pl_\nu v_1||_{N_{m+1}}\le C\delta.$$
	Set
$$\mathcal U_\delta=\Big\{(u_0,u_1,u_2,h)\in \H_m\times A_m\ \Big| (u_0,u_1,u_2,h)~{\rm satisfies}~(\ref{c2})\Big\}.$$
	Taking $\hat w\in Z_\delta:=\{\hat w\in E^{m+2}:\ ||\hat w||_{E^{m+2}}\le \delta\},$
	we define
	$$\mathcal F:\ Z_\delta\to Z_\delta,\quad \mathcal F\hat w=v_2,$$
	where $v_2$ is the solution to the following system:
	\begin{equation}\label{sysv2}
		\left\{\begin{array}{ll}
			v_{2ttt}+\alpha v_{2tt}-b\Delta v_{2t}-c^2\Delta v_2=[p(v_1+\hat w)^2]_{tt},\ & (x,t)\in Q,\\
			v_2=0,\  & (x,t)\in \Sigma,\\
			v_2(x,0)=v_{2t}(x,0)=v_{2tt}(x,0)=0,\  &  x\in \Omega,
		\end{array}
		\right.
	\end{equation}
	and $v_1$ is the unique solution to (\ref{sysv1}).
	Since $p\in E^{m+2}, v_1\in E^{m+2},$  then
	$$[p(x,t)(v_1+\hat w)^2]_{tt}\in E^m.$$

	By the estimate (\ref{estlinear}), and noting that $E^m$ is a Banach algebra for $m>n+1,$ we have
	\beq
	||\mathcal F \hat w||_{E^{m+2}}=||v_2||_{E^{m+2}}&&\le C||[p(v_1+\hat w)^2]_{tt}||_{E^m}\le C||p(v_1+\hat w)^2||_{E^{m+2}}\nonumber\\
	&&\le C||p||_{E^{m+2}}(||v_1||^2_{E^{m+2}}+||\hat w||^2_{E^{m+2}})\le C\delta^2.\nonumber
	\eeq
	Taking $\delta>0$ sufficiently small such that $\delta\le \frac{1}{C}.$ Then $\mathcal F$ is well-defined.
	
	We proceed to prove that $\mathcal F$ is a contraction.
	
	Indeed, for any $\hat w_1,\hat w_2\in Z_\delta$, let $v_2^1=\mathcal F\hat w_1,\ v_2^2=\mathcal F\hat w_2$  respectively solve
	system (\ref{sysv2}) with respect to $(v_1,\hat w_1)$ and $(v_1,\hat w_2).$ Set $v=v_2^1-v_2^2.$
	Then, $v$ solves
	\begin{equation*}
		\left\{\begin{array}{ll}
			v_{ttt}+\alpha v_{tt}-b\Delta v_{t}-c^2\Delta v=[p(v_1+\hat w_1)^2-p(v_1+\hat w_2)^2]_{tt},\ & (x,t)\in Q,\\
			v=0,\  & (x,t)\in \Sigma,\\
			v(x,0)=v_{t}(x,0)=v_{tt}(x,0)=0,\ & x\in \Omega.
		\end{array}
		\right.
	\end{equation*}
	Again using the estimate (\ref{estlinear}), we obtain
	\beq
	||\mathcal F\hat w_1-\mathcal F\hat w_2||_{E^{m+2}}&&\le C||[p(v_1+\hat w_1)^2-p(v_1+\hat w_2)^2]_{tt}||_{E^m}\nonumber\\
	&&\le C||p||_{E^{m+2}}||(\hat w_1+\hat w_2+2v_1)||_{E^{m+2}}||\hat w_1-\hat w_2||_{E^{m+2}}\nonumber\\
	&&\le C\delta ||\hat w_1-\hat w_2||_{E^{m+2}}.\nonumber
	\eeq
	Take $\delta>0$ small enough such that $C\delta<1.$ Then, $\mathcal F:Z_\delta\to Z_\delta$ is a contraction. By the Banach fixed point theorem, we know that (\ref{sysv2}) admits a unique solution $v_2\in E^{m+2},$ which implies that $u=v_1+v_2$ is the unique solution to (\ref{sysmgt}), and (\ref{estnon}) holds.  \qedsymbol
	
	\section{Construction of CGO solutions}
In this section, we construct the CGO solutions for linear M-G-T equation.
In the rest of this paper, and for the convenience of notations, we proceed the calculations and arguments by viewing $(\Omega,g)$ as a Riemannian manifold, with the simple metric $g=b^{-1}ds^2.$
	\subsection{Construction of the phase and amplitude functions}
	Recalling that
	$$\PP u=u_{ttt}+\alpha u_{tt}-b\Delta u_t-c^2\Delta u.$$
	We construct the CGO solutions to $\PP u=0$ of the following form:
	$$u=e^{\boldsymbol{i}\sigma(\var(x)+t)}(a_1(x,t)+\sigma^{-1}a_2(x,t))+R_\sigma(x,t),$$
	where $\boldsymbol{i}$ is the imaginary unit and $a_1,a_2,\var$ are smooth functions, $R_\sigma$ denotes the remainder term, which will be specified and  discussed later.
	
	For any $\psi_1,\psi_2\in H^1(\Omega)$ and $\var_0\in H^2(\Omega),$ we have
	$$\nabla_g\psi_1=b\nabla \psi_1,\ \left\langle\nabla_g\psi_1,\nabla_g\psi_2 \right\rangle_g=b\left\langle\nabla\psi_1,\nabla\psi_2 \right\rangle,\ \Delta_g\var_0=b\Delta\var_0,$$
	where $\Delta_g$ denotes the Beltrami-Laplacian in the metric $g.$
	Directly calculating yields
	\beq\label{solveP}
	&&\PP[e^{\boldsymbol{i}\sigma(\var+t)}(a_1+\sigma^{-1}a_2)]\nonumber\\
	&&=e^{\boldsymbol{i}\sigma(\var+t)}\{\boldsymbol{i}\sigma^3a_1(b|\nabla\var|^2-1)+\sigma^2[(b|\nabla\var|^2-3)a_{1t}+(c^2|\nabla\var|^2-\alpha +b\Del\var)a_1\nonumber\\
	&&\quad +2b\<\nabla\var,\nabla a_1\>+\boldsymbol{i} a_2(|\nabla\var|^2-1) ]+\sigma[\boldsymbol{i}(3a_{1tt}+2\alpha a_{1t}-2\<\nabla\var,c^2\nabla a_1+b\nabla a_{1t}\>\nonumber\\
	&&\quad -(c^2a_1+ba_{1t})\Del\var-b\Del a_1)+(b|\nabla\var|^2-3)a_{2t}+(c^2|\nabla\var|^2-\alpha +b\Del\var)a_2\nonumber\\
	&&\quad+2b\<\nabla\var,\nabla a_2\>]+\boldsymbol{i}[3a_{2tt}+2\alpha a_{2t}-2\<\nabla\var,c^2\nabla a_2+b\nabla a_{2t}\>\nonumber\\
	&&\quad -(c^2a_2+ba_{2t})\Del\var-b\Del a_2]+\PP a_1+\sigma^{-1}\PP a_2 \}.
	\eeq
	
	Collecting the powers of $\sigma^2$ and $\sigma,$ recalling that $\beta=\frac{c^2}{b},$
	we let
	\be\label{trans1}
	(|\nabla_g\var|_g^2-3)a_{1t}+(\beta|\nabla_g\var|_g^2-\alpha +\Del_g\var)a_1+2\<\nabla_g\var,\nabla_ga_1\>_g+\boldsymbol{i} a_2(|\nabla_g\var|_g^2-1)=0,
	\ee
	and
	\beq\label{trans2}
	&&\boldsymbol{i}[3a_{1tt}+2\alpha a_{1t}-2\<\nabla_g\var,\beta\nabla_ga_1+\nabla_g a_{1t}\>_g-(\beta a_1+a_{1t})\Del_g\var-\Del_g a_1]\nonumber\\
	&&\quad+(|\nabla_g\var|_g^2-3)a_{2t}+(\beta|\nabla_g\var|_g^2-\alpha +\Del_g\var)a_2+2\<\nabla_g\var,\nabla_g a_2\>_g=0.
	\eeq
	We choose $\varphi$ so that it satisfies the following eikonal equation:
	\be\label{eikonal}
	|\nabla_g\var|_g=b|\nabla\var|=1.
	\ee
	Clearly, $\varphi(x)=\frac{1}{b}|x-y|$ solves (\ref{eikonal}), where $y\in \mathbb R^n\backslash\overline{\Omega},$ $x\in \Omega.$ However, we calculate and analyze these equations in a general simple Riemannian manifold, discussed e.g., in \cite[Section 2]{YK1}. More precisely,  we firstly extend $(\Omega,g)$ into a simple manifold $(\Omega_1,g)$ in such a
	way that $\Omega$ is contained into the interior of $\Omega_1$.
	
	Let $q\in\partial\Omega_1$ and $\Omega_{1q}$ be the tangent space at $q.$ Considering the polar normal coordinates $(r,\theta)$ on $\Omega_1$ given by
	$$x=\exp_q(r\theta),$$
	where $\exp_q(\cdot)$ is the exponential map on $\Omega_1$ and the parameter $r>0$ denotes the length of the curve connecting $x$ and $q,$
	$$\theta\in S_q\Omega_1=\{\theta\in \Omega_{1q}:\ |\theta|_g=1\}.$$
	It follows from the Gauss lemma that, in
	these coordinates the metric is of the form
	$$g(r,\theta)=dr^2+g_0(r,\theta),$$
	where $g_0(r,\theta)$ is a metric on $S_q\Omega_1$ which depends smoothly on $r.$
	
	In the above setting, we choose
	$$\varphi(x)={\rm dist}_q(x),\quad x\in \Omega,$$
	where ${\rm dist}_q(\cdot)$ denotes the Riemannian distance function from $q$ to $x$ on $(\Omega_1,g).$
	Clearly, $\varphi$ is given by $r$ in the polar normal coordinates, then $\nabla_g\varphi=\frac{\partial}{\partial r}$ and $|\nabla_g\varphi|_g=1.$
	
	We write $a_j(x,t)=a_j(r,\theta,t)$ for $j=1,2,$ and also
	use this notation to indicate the representation in the polar normal coordinates for
	other functions. Notice that in the new coordinate system,
	$$\Delta_g\varphi=\frac{d_r}{2d},$$
	where $d=d(r,\theta)=\det g_0(r,\theta).$

	Recall that $\gamma=\alpha-\frac{c^2}{b}.$ Then, (\ref{trans1}) and (\ref{trans2}) can be reduced to the following two transport equations:
	\be\label{trans3}
	a_{1t}-a_{1r}+\Big(\frac{\gamma}{2}-\frac{d_r}{4d}\Big)a_1=0,
	\ee
	\be\label{trans4}
	a_{2t}-a_{2r}+\Big(\frac{\gamma}{2}-\frac{d_r}{4d}\Big)a_2=\tilde{\xi},
	\ee
	where
	$$2\tilde{\xi}(x,t)=\boldsymbol{i}[3a_{1tt}+2\alpha a_{1t}-2\<\nabla_g\var,\beta\nabla_ga_1+\nabla_g a_{1t}\>_g-(\beta a_1+a_{1t})\Del_g\var-\Del_g a_1].$$
	For any $\chi\in C^{\infty}(\mathbb R),$ and $\phi\in C^{\infty}(S_y\Omega_1),$ it is easy to see
	\be\label{solvea1}
	a_1(r,\theta,t)=e^{-\frac\mu 2(r+t)}\chi(r+t)\phi(\theta)e^{-\frac\gamma 2 t}d^{-\frac14}
	\ee
	solves (\ref{trans3}), where $\mu>0$ is a constant. It is not necessary to give the explicit expression of $a_2$ when carrying on our analysis, since the factor $\sigma^{-1}a_2$ vanishes as $\sigma\to\infty.$
	\begin{rem}
		Generally, the non-homogeneous transport equation {\rm(\ref{trans4})}, even for the higher dimensional case, is globally solvable, see {\rm\cite[Proposition 4]{RDJL}}. Indeed,
		\beq\label{transhom}
		a_2(r,\theta,t)=\psi(r+t)\phi(\theta)e^{-\int_0^t\zeta(r+s,\theta)ds}+\int_0^t\tilde{\xi}(r+t-s,\theta,s)e^{-\int_s^t\zeta(r+l,\theta)dl}ds,
		\eeq
		where $\psi\in C^{\infty}(\mathbb R)$ is arbitrary and $\zeta(r,\theta)=(\frac\gamma 2-\frac{d_r}{4d}).$
		Moreover,  the a prior estimate of {\rm(\ref{trans4})} can be found in {\rm\cite[Proposition 2.4]{LUH}}. That is, there exist two positive constants $C_1,C_2$ such that
		\begin{eqnarray}\label{a2bouna1}
			\begin{array}{ll}
				||a_2||_{L^\infty(Q)}\le C_1(M,T,a_2(x,0))+C_2(M,T)||\tilde{\xi}||_{L^\infty}.
			\end{array}
		\end{eqnarray}
	\end{rem}
	\subsection{Asymptotic analysis of the remainder term}
	By the above discussions on the construction of the CGO solutions to $\PP u=0,$ we are left to deal with the remainder term $R_\sigma$ satisfying
	\begin{equation}\label{remeqn}
	\left\{\begin{array}{ll}
		R_{\sigma ttt}+\alpha R_{\sigma tt}-b\Delta R_{\sigma t}-c^2\Delta R_{\sigma}=-e^{\boldsymbol{i}\sigma t}F_\sigma(a_1,a_2,\varphi),\ & (x,t)\in Q,\\
		R_\sigma=0,\ & (x,t)\in \Sigma,\\
		R_{\sigma}(x,0)=R_{\sigma t}(x,0)=R_{\sigma tt}(x,0)=0,\ & x\in \Omega,
	\end{array}
	\right.
	\end{equation}
	where
	\beq\label{Fsigma}
	F_\sigma(a_1,a_2,\varphi)&&=e^{\boldsymbol{i}\sigma\var}\{\boldsymbol{i}[3a_{2tt}+2\alpha a_{2t}-2\left\langle\nabla_g\var,\beta\nabla_g a_2+\nabla_ga_{2t}\right\rangle_g-(\beta a_2+a_{2t})\Delta_g\var]\nonumber\\
	&&\quad +\PP a_1+\sigma^{-1}\PP a_2\}\in H^1(0,T;L^2(\Om)).
	\eeq
	
	Now, we are in a position to state the following theorem, which is useful in analyzing the inverse problems.
	\begin{thm}\label{construccgo}
		The equation $\PP u=0$ admits a CGO solution of the following form
		$$u=e^{\boldsymbol{i}(\var+t)}(a_1+\sigma^{-1}a_2)+R_\sigma\in E^2,$$
		where $R_\sigma$ satisfies {\rm(\ref{remeqn})}. Moreover, there exists a positive constant $C=C(T,M,\Omega)$ such that
		\beq\label{asyrem}
		&&\sigma\left(||R_\sigma||_{L^2(Q)}+||R_{\sigma t}||_{L^2(Q)}+||\nabla R_\sigma||_{L^2(Q)}\right) +||R_{\sigma tt}||_{L^2(Q)}+||\nabla R_{\sigma t}||_{L^2(Q)}\nonumber\\
		&&\le C||F_\sigma(a_1,a_2,\var)||_{H^1(0,T;L^2(\Om))},
		\eeq
		which implies that $||R_\sigma||_{L^2(Q)}$ and $||R_{\sigma t}||_{L^2(Q)}$ vanish as $\sigma\to \infty.$
	\end{thm}
	
	{\bf Proof.}\quad The well-posedness result follows from Theorem \ref{well-lower}. Recalling the expression of $F_\sigma(a_1,a_2,\varphi),$ applying (\ref{estlinear}) in Theorem \ref{well-linear} with $R_\sigma=v,$ we have
	$$||\nabla R_{\sigma t}||_{L^2(Q)}+||R_{\sigma tt}||_{L^2(Q)}\le C||e^{\boldsymbol{i}\sigma t}F_\sigma(a_1,a_2,\varphi)||_{L^2(Q)}\le C||F_\sigma(a_1,a_2,\varphi)||_{L^2(Q)}.$$
	It remains to show
	$$\sigma\left(||R_{\sigma t}||_{L^2(Q)}+||\nabla R_\sigma||_{L^2(Q)}+||R_\sigma||_{L^2(Q)}\right)\le C||F_\sigma(a_1,a_2,\varphi)||_{H^1(0,T;L^2(\Omega))}.$$
	Indeed, let $v_\sigma(x,t)=\int_0^tR_\sigma(x,l)dl.$ Then, $v_\sigma$ solves
	\begin{eqnarray}\label{sysvvv}
		\left\{ \begin{array}{ll}
			v_{\sigma ttt}+\alpha v_{\sigma tt}-b\Del v_{\sigma t}-c^2\Del v_{\sigma}=\mathcal L_\sigma(x,t) & \mbox{ in } Q,\\
			v_\sigma=0  & \mbox{ on } \Sigma,\\
			v_\sigma(x,0)=v_{\sigma t}(x,0)=v_{\sigma tt}(x,0)=0 & \mbox{ in } \Om,
		\end{array} \right.
	\end{eqnarray}
	where
	$$\mathcal L_\sigma(x,t)=-\int_0^te^{\boldsymbol{i}\sigma l}F_\sigma(a_1,a_2,\varphi)(x,l)dl.$$
	Respectively multiplying the first equation of (\ref{sysvvv}) by $v_{\sigma tt}$ and $v_{\sigma t},$ and integrating by parts over $\Om\times(0,t),$ respectively yields
	\beq\label{ener1}
	&&\dis\frac12\int_\Om(v_{\sigma tt}^2+|\nabla_g v_{\sigma t}|_g^2)dx+\alpha\int_0^t\int_\Omega|v_{\sigma tt}|^2dxdl-\beta\int_0^t\int_\Om|\nabla_g v_{\sigma t}(x,l)|_g^2dxdl\nonumber\\
	&&\quad +\beta\int_\Om\<\nabla_g v_{\sigma t},\nabla_g v_\sigma\>_gdx=\int_0^t\int_\Om\mathcal L_\sigma(x,l)v_{\sigma tt}(x,l)dxdl,
	\eeq
	and
	\beq\label{ener2}
	&&\dis\frac12\int_\Om(\alpha v_{\sigma t}^2+\beta|\nabla_g v_\sigma|_g^2)dx+\int_0^t\int_\Om|\nabla_g v_{\sigma t}(x,l)|_g^2dxdl-\int_0^t\int_\Om v_{\sigma tt}^2(x,l)dxdl\nonumber\\
	&&\quad +\int_\Om v_{\sigma t}v_{\sigma tt}dx=\int_0^t\int_\Om\mathcal L_\sigma(x,l)v_{\sigma t}(x,l)dxdl.
	\eeq
	Combining (\ref{ener1}) with (\ref{ener2}) and noting that
	$$
	v_{\sigma t}(x,t)=\int_0^tv_{\sigma tt}(x,l)dl,\quad \nabla v_\sigma(x,t)=\int_0^t\nabla v_{\sigma t}(x,l)dl,
	$$
	it holds that
	\begin{eqnarray*}\label{ener0}
	\hat E(t)&&\dis:=\frac12\int_\Om(v_{\sigma tt}^2+|\nabla_g v_{\sigma t}|_g^2+|\nabla_g v_\sigma|_g^2+v_{\sigma t}^2)dx\nonumber\\
	&&\dis\le C\int_0^t\hat E(l)dl+\int_0^t\int_\Om\mathcal L_\sigma(x,l)(v_{\sigma t}+v_{\sigma tt})(x,l)dxdl.
	\end{eqnarray*}
	By some similar arguments with \cite[Lemma 3.1]{HGH} (see also \cite[Lemma 2.2]{AFYK}), we have
	\begin{eqnarray*}\label{ener3}
	&&\int_0^t\int_\Om\mathcal L_\sigma(x,l)(v_{\sigma t}+v_{\sigma tt})(x,l)dxdl\nonumber\\
	&&=-\int_0^t\int_\Om\[\int_0^le^{\boldsymbol{i}\sigma s}F_\sigma(a_1,a_2,\varphi)(x,s)ds \](v_{\sigma t}+v_{\sigma tt})(x,l)dxdl\nonumber\\
	&&\le C\|\|\int_0^te^{\boldsymbol{i}\sigma l}F_{\sigma}(a_1,a_2,\varphi)(x,l)dl\|\|^2_{L^2(Q)}+C(T)\int_0^t \hat E(l)dl.
	\end{eqnarray*}
	Therefore, we obtain
	\begin{eqnarray*}
		\begin{array}{ll}
			\dis \hat E(t)\le C \int_0^t \hat E(l)dl+C\|\|\int_0^te^{\boldsymbol{i}\sigma l}F_\sigma(a_1,a_2,\varphi)(x,l)dl\|\|^2_{L^2(Q)}.
		\end{array}
	\end{eqnarray*}
	Invoking the Gronwall inequality yields
	\begin{eqnarray*}
		\begin{array}{ll}
			\dis \hat E(t)\le C\|\|\int_0^te^{\boldsymbol{i}\sigma l}F_{\sigma}(a_1,a_2,\varphi)(x,l)dl\|\|^2_{L^2(Q)}.
		\end{array}
	\end{eqnarray*}
	Recalling that $v_{\sigma t}=R_\sigma, v_{\sigma tt}=R_{\sigma t},$ we have
	\beq\label{sigmare}
	&&||R_{\sigma t}||^2_{L^2(Q)}+||\nabla R_\sigma||^2_{L^2(Q)}+||R_\sigma||^2_{L^2(Q)}\nonumber\\
	&&\le C\int_Q\|\int_0^te^{\boldsymbol{i}\sigma l}F_\sigma(a_1,a_2,\var)(x,l)dl\|^2dxdt\nonumber\\
	&&=\frac{C}{\sigma^2}\!\int_Q\|\int_0^t(\pl_le^{\boldsymbol{i}\sigma l})F_\sigma(a_1,a_2,\var)(x,l)dl\|^2dxdt\nonumber\\
	&&\le \frac{C}{\sigma^2}\int_Q|e^{\boldsymbol{i}\sigma t}F_\sigma(a_1,a_2,\var)(x,t)-F_\sigma(a_1,a_2,\var)(x,0)|^2dxdt\nonumber\\
	&&\quad+\frac{C}{\sigma^2}\int_Q\|\int_0^t e^{\boldsymbol{i}\sigma l}\pl_lF_\sigma(a_1,a_2,\var)(x,l)dl\|^2dxdt\nonumber\\
	&&\le \frac{C}{\sigma^2}\int_Q\left[ \|e^{\boldsymbol{i}\sigma t}F_\sigma(a_1,a_2,\var)(x,t)\|^2+\|\int_0^t e^{\boldsymbol{i}\sigma l}\pl_lF_\sigma(a_1,a_2,\var)(x,l)dl\|^2\right] dxdt,\nonumber
	\eeq
	where in the last inequality,
	$$F_\sigma(a_1,a_2,\varphi)(x,0)=F_\sigma(a_1,a_2,\varphi)(x,t)-\int_0^t\partial_lF_\sigma(a_1,a_2,\varphi)(x,l)dl$$
	is used.
	In view of $F_\sigma(a_1,a_2,\varphi)$ given by (\ref{Fsigma}), the expression (\ref{solvea1}) and the a prior estimate of $a_2$ bounded by $a_1$ (see (\ref{a2bouna1})), notice that $F_\sigma(a_1,a_2,\varphi)(x,t),\ \partial_tF_\sigma(a_1,a_2,\varphi)(x,t)$
	are uniformly continuous with respect to time $t\in [0,T],$ and thus bounded over $Q,$
	we have
	$$ ||R_{\sigma t}||^2_{L^2(Q)}+||\nabla R_\sigma||^2_{L^2(Q)}+||R_\sigma||^2_{L^2(Q)}\le \frac{C}{\sigma^2}||F_\sigma(a_1,a_2,\varphi)||^2_{H^1(0,T;L^2(\Omega))}.$$
	Therefore, the proof of estimate (\ref{asyrem}) is completed. \qedsymbol

	\section{Proofs of the main results}
	This section concerns to prove the main results of this present paper. The second order linearization below, and the CGO solutions constructed in Section 4 are used to derive the uniqueness of recovering the nonlinearity $p(x,t).$
	
	Let us introduce below the second order linearization procedure to the nonlinear system (\ref{sysmgt}).
	\subsection{Second order lineartization}
	Let $\varepsilon=(\varepsilon_1,\varepsilon_2)\in \mathbb N^2$ and let 
	$$h_j\in A_m,\ (u_{0j},u_{1j},u_{2j})\in H^{m+2}(\Omega)\times H^{m+1}(\Omega)\times H^m(\Omega),\ m>n+1, j=1,2.$$
	Assume that $u(x,t;\varepsilon)$ is the solution to the nonlinear system (\ref{sysmgt}) with respect to
	\begin{eqnarray*}
		\begin{array}{ll}
			{ Boundary\ data}:~ h=\varepsilon_1h_1+\varepsilon_2h_2,\quad (x,t)\in\Sigma,&\\[2mm]
			{ Initial\ data}:~ u_0=\varepsilon_1u_{01}+\varepsilon_2u_{02},\ u_1=\varepsilon_1u_{11}+\varepsilon_2u_{12},\ u_2=\varepsilon_1u_{21}+\varepsilon_2u_{22},\ x\in\Omega.&
		\end{array}
	\end{eqnarray*}
	The parameter $|\varepsilon|=|(\varepsilon_1,\varepsilon_2)|$ is sufficiently small such that
	$$||h||_{A_m}+||u_0||_{H^{m+2}(\Omega)}+||u_1||_{H^{m+1}(\Omega)}+||u_2||_{H^m(\Omega)}<\delta.$$
	Then, it follows from Theorem \ref{wellnonlinear} that $u\in E^{m+2}$ and  $u(x,t;0)=0.$

\ms
	
	Let $w_j=\frac{\partial u(x,t;\varepsilon)}{\partial\varepsilon_j}\Big|_{\varepsilon=0}$ for $j=1,2.$ Then for each $j=1,2,$ $w_j$ solves
	\be\label{linearcgo}
	\left\{
	\begin{array}{ll}
		w_{jttt}+\alpha w_{jtt}-b\Del w_{jt}-c^2\Del w_j=0,\  &  (x,t)\in Q,\\
		w_j=h_j,\ & (x,t)\in \Sigma,\\
		w_j(x,0)=u_{0j}(x),w_{jt}(x,0)=u_{1j}(x),w_{jtt}(x,0)=u_{2j}(x),\ &  x\in \Om.
	\end{array}
	\right.
	\ee
	Let $w=\frac{\partial^2u(x,t;\varepsilon)}{\partial\varepsilon_1\partial\varepsilon_2}\Big|_{\varepsilon=0}.$ Then $w$ solves
	\begin{equation*}
	\left\{
	\begin{array}{ll}
		w_{ttt}+\alpha w_{tt}-b\Del w_{t}-c^2\Del w=F(p,w_1,w_2),\ & (x,t)\in Q,\\
		w=0,\ & (x,t)\in \Sigma,\\
		w(x,0)=w_{t}(x,0)=w_{tt}(x,0)=0,\  & x\in \Om,
	\end{array}
	\right.
	\end{equation*}
	where
	\beq\label{linearu12}
	F(p,w_1,w_2)&&=2p(w_{1tt}w_2+2w_{1t}w_{2t}+w_1w_{2tt})\nonumber\\
	&&\quad+4p_t(w_{1t}w_2+w_1w_{2t})+2p_{tt}w_1w_2.\nonumber
	\eeq
	\subsection{Proof of Theorem \ref{thm1}}
	Based on the second order linearization, together with the construction of CGO solutions to $Pw_j$=0 for $j=1,2,$ we are ready to prove the first result.

\ms
	
	{\bf Proof of Theorem \ref{thm1}}\quad
	Recall that $\gamma=\alpha-\frac{c^2}{b}=\alpha-\beta.$ Let $W=w_t+\beta w.$ Then, $W$ is the solution to the following system:
	\be\label{WWW}
	\left\{
	\begin{array}{ll}
		W_{tt}-b\Del W+\gamma W_t=F(p,w_1,w_2)+\gamma\beta w_t,\ & (x,t)\in Q,\\
		W=0,\ & (x,t)\in \Sigma,\\
		W(x,0)=W_{t}(x,0)=0,\ & x\in \Om.
	\end{array}
	\right.
	\ee
	Let $y\in C([0,T];H^1(\Omega))\cap C^1([0,T];L^2(\Omega))\subset L^2(Q)$ be a solution to
	\be\label{sysy}
	L_0y:=y_{tt}-b\Del y-\gamma y_t=0,\ (x,t)\in Q.
	\ee
	
	Multiplying the first equation of (\ref{WWW}) by $y$ and integrating over $\Omega\times (0,T),$ applying the divergence theorem, we have
	\beq\label{muly}
	&&\int_0^T\int_\Om[F(p,w_1,w_2)+\gamma\beta w_t]ydxdt\nonumber\\
	&&=b\int_0^T\int_\Ga y\pl_\nu Wd\Ga dt+\int_\Om(yW_t-y_tW+\gamma yW)dx\Big|_{t=T}.
	\eeq
	To prove that the measurement map $\Lambda_T$ uniquely determines $p(x,t),$ it suffices to prove $\Lambda_T(h,u_0,u_1,u_2)=0$ implies that $p(x,t)=0$ in $Q.$
	Then, by the assumption that
	$$w(x,T)=w_t(x,T)=w_{tt}(x,T)=0\quad {\rm for}\ x\in \Omega,\quad \partial_\nu w=0\quad {\rm for}\ (x,t)\in \Sigma,$$
	which gives
	$$W(x,T)=W_t(x,T)=0\quad{\rm for}\ x\in \Omega,\quad \partial_\nu W=0\quad {\rm for}\ (x,t)\in \Sigma.$$
	From (\ref{muly}), we have,
	\beq\label{intiden}
	\int_0^T\int_\Omega yF(p,w_1,w_2)dxdt=-\gamma\beta\int_0^T\int_\Om yw_tdxdt.
	\eeq
	
	By the construction of CGO solutions in Section 4, with $\sigma$ respectively replacing by $2\sigma$ and $-\sigma$ in the expression (\ref{solveP}),  applying Theorem \ref{construccgo}, we respectively choose
	$$w_1=e^{2\boldsymbol{i}\sigma (\varphi+t)}(a_{11}+(2\sigma)^{-1}a_{12})+R_{1\sigma},$$
	$$w_2=e^{-\boldsymbol{i}\sigma (\varphi+t)}(a_{21}+\sigma^{-1}a_{22})+R_{2\sigma},$$
	is the CGO solutions to the linear system (\ref{linearcgo}) for $j=1,2.$
	To the above CGO solutions, we choose
	$$a_{11}(r,\theta,t)=e^{-\frac{\mu(r+t)}{2}}\phi(\theta)e^{-\frac\gamma 2 t}d(r,\theta)^{-\frac14},\quad a_{21}(r,\theta,t)=e^{-\frac{\mu(r+t)}{2}}e^{-\frac\gamma 2 t}d(r,\theta)^{-\frac14}.$$

	By some similar arguments as the construction of CGO solutions to $\PP u=0$ proposed in Section 4, it is easy to see that
	$$y=e^{-\boldsymbol{i}\sigma(\var+t)}a_0+r_{0\sigma}$$
	is a solution to (\ref{sysy}), and $r_{0\sigma}$ solves
	\begin{equation*}
	\left\{\begin{array}{ll}
		r_{0\sigma tt}-b\Del r_{0\sigma}-\gamma r_{0\sigma t}=-e^{-\boldsymbol{i}\sigma t}Y_\sigma(a_0,\var),\ & (x,t)\in Q,\\
		r_{0\sigma}=0,\  & (x,t)\in \Sigma,\\
		r_{0\sigma}(x,0)=r_{0\sigma t}(x,0)=0,\ & x\in \Omega,
	\end{array}
	\right.
	\end{equation*}
	where
	$$Y_\sigma(a_0,\varphi)=e^{-\boldsymbol{i}\sigma\varphi}(a_{0tt}-b\Delta a_0-\gamma a_{0t})\in L^2(Q).$$
	Moreover, similar to the proof of Theorem \ref{construccgo} (see also \cite[Lemma 2.2]{AFYK} or \cite[Lemma 2.2]{YK1}), we can obtain
	\be\label{r0}
	\lim_{\sigma\to\infty}||r_{0\sigma}||_{L^2(Q)}=0.
	\ee

	Inserting the CGO solutions into the expression $F(p,w_1,w_2),$ yields
	$$F(p,w_1,w_2)=e^{\boldsymbol{i}\sigma t}[-\sigma^2a_{11}a_{12}+O(\sigma^2)(R_{1\sigma}+R_{2\sigma})]e^{\boldsymbol{i}\sigma\varphi}+\mathcal O(\sigma).$$
	Therefore, by some similar arguments in obtaining the estimate (\ref{asyrem}), we have
	\be\label{asyw}
	\lim_{\sigma\to\infty}\sigma^{-2}||w_t||_{L^2(Q)}=0.
	\ee
	Let $\sigma\to\infty.$ Notice that (\ref{asyrem}) for $R_{1\sigma},R_{2\sigma},$ (\ref{r0}) for $r_{0\sigma},$ and (\ref{asyw}) for $w_t.$ It then follows from (\ref{intiden}) that
	\beq\label{a1a2}
	\int_0^T\int_\Om a_{11}a_{12}(pa_0)dxdt=-\gamma\beta\lim_{\sigma\to\infty}\sigma^{-2}\int_0^T\int_\Om yw_tdxdt=0.
	\eeq
	
	We extend $p$ by zero to $\Omega_1\times(0,+\infty)$ and still denoted by $p$. Denote
	$$\tau_+(q,\theta)=\inf\{r>0:\gamma_{q,\theta}(r)\in \partial\Omega_1\},$$
	where $\gamma_{q,\theta}$ is the geodesic on $\Omega_1$ initiating from $q$ with $\dot{\gamma}_{q,\theta}(0)=\theta.$ From (\ref{a1a2}), and noting that
	$dx=d(r,\theta)^{\frac12}drd\theta,$ we obtain in the polar normal coordinates that
	$$
	\int_0^{+\infty}\int_{S_q\Omega_1}\int_0^{\tau_+(q,\theta)}\t p(r,\theta,t)\phi(\theta)e^{-\mu(r+t)}drd\theta dt=0,
	$$
	for all $\phi\in C^\infty(S_q\Omega_1),$ $q\in\partial\Omega_1$ and $\mu>0,$ where $$\t p(r,\theta,t)=e^{-\gamma t}(pa_0)(r,\theta,t)\in L^1(0,+\infty)$$ with respect to $t.$ By the same arguments as that in \cite[Section 2]{YK1}, we conclude that
	$pa_0=0$ in $Q.$ By taking
	$$a_0=e^{\frac{\gamma t}{2}}d(r,\theta)^{-\frac{1}{4}}\ne 0,$$
	we obtain that $p=0.$ Thus, the proof of Theorem \ref{thm1} is finished. \qedsymbol
	
	\subsection{Proof of Theorem \ref{thm2}}
	Let the boundary data $u(x,t)|_{\Sigma}=h(x,t)=\varepsilon_1h_1(x,t)+\varepsilon_2h_2(x,t)$ be the system input only.
	For each $j=1,2,$ $w_j=\frac{\partial u(x,t;\varepsilon)}{\partial\varepsilon_j}\Big|_{\varepsilon=0}$ satisfies
	\be\label{linewj}
	\left\{
	\begin{array}{ll}
		w_{jttt}+\alpha w_{jtt}-b\Del w_{jt}-c^2\Del w_j=0,\ & (x,t)\in Q,\\
		w_j=h_j,\ & (x,t)\in \Sigma,\\
		w_j(x,0)=w_{jt}(x,0)=w_{jtt}(x,0)=0,\ & x\in \Om.
	\end{array}
	\right.
	\ee
	
The CGO solutions can not be directly used due to the zero initial data of the above system (\ref{linewj}). To avoid this, we suitably choose the CGO solutions to (\ref{linewj}) such that they vanish at $t=0.$
Indeed, we choose
$$a_{11}=e^{-\frac{\mu}{2}(r+t)}\chi(r+t)\phi(\theta)e^{-\frac\gamma 2 t}d(r,\theta)^{-\frac14},\quad a_{21}=e^{-\frac{\mu}{2}(r+t)}\chi(r+t)e^{-\frac\gamma 2 t}d(r,\theta)^{-\frac14},$$
where $\chi\in C^{\infty}(0,+\infty)$ is a cut-off function satisfying
\begin{equation*}
\chi(s)=\left\{\begin{array}{ll}
0,\ & s\in (0,T^*+\t\delta],\\
0< \chi(s)\le 1,\ & s\in (T^*+\t\delta,T^*+2\t\delta),\\
1,\ & s\in [T^*+2\t\delta,+\infty).
\end{array}
\right.
\end{equation*}
Notice that $r\le T^*+\t\delta,$ then it is clear that $$\partial_t^ka_{11}(r,\theta,t)|_{t=0}=\partial_t^ka_{21}(r,\theta,t)|_{t=0},\quad k=0,1,2.$$
Similar to the above $a_{11}$ and $a_{21},$ using the equations (\ref{trans4})  and expressions (\ref{transhom}) for $a_{12}, a_{22}$, and choosing $\psi=\chi$, we see
$$\partial_t^ka_{12}(r,\theta,t)|_{t=0}=\partial_t^ka_{22}(r,\theta,t)|_{t=0}=0,\quad k=0,1,2.$$
Together with $\partial_t^kR_{1\sigma}(x,t)|_{t=0}=\partial_t^kR_{2\sigma}(x,t)|_{t=0}=0,$ $k=0,1,2,$ we know that the CGO solutions
$$w_1=e^{2\boldsymbol{i}\sigma(\varphi+t)}(a_{11}+(2\sigma)^{-1}a_{12})+R_{1\sigma},\quad w_2=e^{-\boldsymbol{i}\sigma(\varphi+t)}(a_{21}+(2\sigma)^{-1}a_{22})+R_{2\sigma}$$
satisfy $w_j(x,0)=w_{jt}(x,0)=w_{jtt}(x,0)=0,\ j=1,2.$

\ms

{\bf Proof of Theorem \ref{thm2}}\quad
By a similar proof with Theorem \ref{thm1}, with the same $y,$ and from $(\ref{a1a2})$, we conclude that
\beq
&&0=\int_0^{+\infty}\int_{S_q\Omega_1}\int_0^{\tau_+(q,\theta)}\t p(r,\theta,t)\chi^2(r+t)\phi(\theta)e^{-\mu(r+t)}drd\theta dt\nonumber\\
&&=\int_{T^*+2\t\delta}^{T}\int_{S_q\Omega_1}\int_0^{\tau_+(q,\theta)}\t p(r,\theta,t)\phi(\theta)e^{-\mu(r+t)}drd\theta dt\nonumber\\
&&=\int_0^T\int_{S_q\Omega_1}\int_0^{\tau_+(q,\theta)}\t p(r,\theta,t)\phi(\theta)e^{-\mu(r+t)}drd\theta dt,\nonumber
\eeq
where $\t p(r,\theta,t)=e^{-\gamma t}(pa_0)(r,\theta,t)$, ${\rm supp}\ p\subset\Omega\times(T^*+2\t\delta,T),$ and $\chi(s)=0$ for $s\le T^*+\t\delta$ and $\chi(s)=1$ for $s\in [T^*+2\t\delta,+\infty)$ are used.

Similar to the proof of Theorem \ref{thm1}, we can finally obtain $p=0.$ Thus, the proof of Theorem \ref{thm2} is completed. \qedsymbol

\begin{rem}
{We mention that, in {\rm\cite[Theorem 5.1]{LLL}}, the authors proved a Runge approximation result so that the CGO solutions can be approximated. That is, the set of solutions to a wave equation with zero initial data is dense in the set of the free wave equation with domain $\Omega\times(t_1,t_2),$ where $2{\rm diam}\Omega<t_1<t_2<T.$ Similarly, the same dense result can also be obtained for the M-G-T equation, if the observability inequality is established for the M-G-T equation. However, 
	it seems hard to obtain the desired observability inequality via the classical Carleman estimates (see \cite[Page 4]{RRAS}). We  here do not intend to discuss more in this present paper.
}
\end{rem}

\section{Concluding remarks}

In this paper, we have proved the uniqueness of determining the nonlinearity appearing in the J-M-G-T equation under two kinds of system inputs. The main results are derived by the second order linearization to the nonlinear equation, and construction of CGO solutions for the linearized equation.

There are still many questions worth considering. An interesting work would be the study of the partial data inverse problems of uniquely determining the time-dependent coefficients $\alpha,b,c,$ $p,$ and it is more challenging to recover them  simultaneously by using less measurements. Also, the stability of recovering these parameters for the nonlinear J-M-G-T equation is worth to be considered.

\end{document}